\documentclass[11pt]{article}

\usepackage{amsfonts}
\usepackage{amssymb}
\usepackage{latexsym}
\usepackage[centertags]{amsmath}
\usepackage{amssymb}
\usepackage{color}

\makeatletter\@addtoreset{equation}{section} \makeatother

\renewcommand\thefigure{\thesection.\@arabic\c@figure}
\renewcommand\thetable{\thesection.\@arabic\c@table}

      \newtheorem{theorem}{Theorem}[section]

      \newtheorem{lemma}[theorem]{Lemma}

      \def\nn{\nonumber}
      \def\rf#1{\mbox{$(\ref{#1})$}}

      \def\be{\begin{equation}} 
      \def\ee{\end{equation}} 
      \def\beqn{\begin{eqnarray}} 
      \def\eeqn{\end{eqnarray}} 
      \def\beq{\begin{eqnarray*}} 
      \def\eeq{\end{eqnarray*}}
      \def\proof{{\noindent\bf Proof\quad}\ }
      \def\mb{\mbox} 
      \def\ra{\rightarrow} 

\begin{document}

\title{ Poisson-Dirichlet Distribution with Small Mutation Rate}

\author{Shui Feng\thanks{Research supported by
      the Natural Science and Engineering Research Council of Canada}
      \\Department of Mathematics and Statistics\\ McMaster
      University\\ Hamilton, Ontario, Canada L8S 4K1}

\date{\today}
\maketitle
\begin{abstract}
The behavior of the Poisson-Dirichlet distribution with small mutation rate is studied through large deviations.
The structure of the rate function indicates that the number of alleles is finite
 at the instant when mutation appears. The large deviation results are then used to study the asymptotic behavior of the
 homozygosity, and the Poisson-Dirichlet distribution with symmetric selection. The latter shows that
 several alleles can coexist when selection intensity goes to infinity in a particular way
 as the mutation rate approaches zero.
\end{abstract}

\vspace*{.125in} \noindent {\bf Key words:} Poisson-Dirichlet
      distribution, Dirichlet processes, homozygosity, large deviations, selection.
      \vspace*{.125in}


\section{Introduction}
\setcounter{equation}{0}

For $\theta >0$, let
$V_1(\theta)\geq V_2(\theta)\geq\cdots$ be the points of a
nonhomogeneous Poisson process with mean measure density
\[
\theta v^{-1}e^{-v}, v >0.
\]
Set
$$V(\theta) =\sum_{i=1}^{\infty}V_i(\theta),$$ and
\be\label{intro1} {\bf P}(\theta )=(P_1(\theta),
P_2(\theta),...)=\left(\frac{V_1(\theta)}{V(\theta)},\frac{V_2(\theta)}{V(\theta)},\cdots\right).
\ee Then ${\bf P}(\theta )$ and $V(\theta)$
are independent, and $V(\theta)$ is a $Gamma(\theta,1)$-distributed random
variable. The law of ${\bf P}(\theta)$ is called the Poisson-Dirichlet
distribution with parameter $\theta$, and is denoted by
$PD(\theta)$.

The labeled version of the Poisson-Dirichlet distribution, called the Dirichlet process, is defined as the law of
\be \label{dirichlet-process-def} \Xi_{\theta,
\nu}=\sum_{k=1}^{\infty}P_k(\theta)\delta_{X_k},\ee
where $X_k, k=1,...$ is a sequence of i.i.d. random variables,
independent of ${\bf P}(\theta )$, with a common diffusive
distribution $\nu$ on $[0,1]$, i.e., $\nu(\{x\})=0$ for every $x$ in
$[0,1]$. The Dirichlet process was introduced
      in \cite{Fer73} as a prior distribution in the context of Bayesian statistics.

    The Poisson-Dirichlet distribution was introduced
      by Kingman \cite{kingman75} to describe the distribution of gene frequencies in a large neutral
      population at a particular locus with each component $P_k(\theta)$ representing the
      proportion of the $k$th most frequent allele. It is the unique reversible measure (cf. \cite{EK81}) of the
      infinitely many neutral alleles diffusion process with state space
      \[
        \nabla =\{{\bf p}=(p_1,p_2,\cdots): p_1\geq p_2\geq\cdots\geq 0, \sum_{i=1}^{\infty}p_i \leq1\},
      \]
      and generator
      \[
      A =\frac{1}{2}\sum_{i,j=1}^{\infty}p_i(\delta_{ij}-p_j)\frac{\partial^2}{\partial p_i \partial p_j}-
      \frac{\theta}{2}\sum_{i=1}^{\infty}p_i\frac{\partial}{\partial p_i}
      \]
      defined on an appropriate domain.

     The parameter $\theta$ represents the scaled mutation rate of a population in the context of population genetics.
     In Bayesian statistics, it can be interpreted as the prior sample size. When $\theta$ is large, the proportions of
     different alleles under $PD(\theta)$ are evenly spread and approach zero; while for small values of $\theta$,
     $PD(\theta)$
     will concentrate mostly around the point $(1,0,\cdots)$. There are extensive studies of the asymptotic
     behavior of $PD(\theta)$ when $\theta$ goes to infinity
     (\cite{watterson-guess77},\cite{griffiths79},\cite{JKK02},\cite{dawson-feng06}, \cite{fenggao08}).
     Since the proportions of alleles are evenly spread and uniformly small, it is thus natural to see
     Gaussian structures (\cite{JKK02},\cite{fenggao08}) for
     large $\theta$. For small mutation rates, the study is very limited. The author is aware of only results
     in \cite{sethti82}
     for $Dirichlet(\theta,\nu)$, and in \cite{eth81} and \cite{ethgri93} for the
     infinitely many neutral alleles diffusion model.

    The case of $\theta=1$ is special. It appears as
     an asymptotic distribution in random number theory (\cite{ABT03}). It is also a critical value in the boundary
     behavior of the infinitely many neutral alleles model.  By using techniques from the theory of
     Dirichlet forms, it was shown in \cite{schmu91} that for the infinitely many neutral alleles model,
     with probability one, there will exist times at which the sample path will hit the boundary of a
     finite-dimensional sub-simplex of $\nabla$ or, equivalently, the single point $(1,0,\cdots)$ iff
     $\theta$ is less than one. The intuition here is that it is possible to have finite number of
     alleles in the population if mutation rate is small.

     But in equilibrium, with $PD(\theta)$ probability one, the number of alleles is always infinity as long
     as $\theta$ is strictly positive. In other words, the critical value between finite number of alleles and
     infinite number of alleles is zero for $PD(\theta)$. In physical terms this sudden change from one to infinity
     can be viewed as a phase transition. The objective of this paper is to investigate
     the microscopic structures during this phase transition. The tool we use is from the theory of large deviations.
     Intuitively, it is unlikely to get a large number of alleles when the mutation first appears.
     Our result will confirm this intuition rigorously, and reveals a transition structure that
     can be viewed as a ``ladder of energy".

  The paper is organized as follows.
 In Section 2, we establish the large deviation principle for $PD(\theta)$ when $\theta$ goes to zero.
The rate function is identified explicitly. In Section 3, the large deviation result is applied to study the asymptotic
behavior of the homozygosity and the impact of selection or exponential tilting. It will be shown that,
in contrast to the neutral case, the population under overdominant selection can preserve more than one alleles
when the mutation rate goes to zero and the selection intensity goes to infinity in a particular way.

\section{Large Deviations}

In this section, we establish the large deviation principle for $PD(\theta)$ when $\theta$ goes to zero.
The result will be obtained through a series of lemmas and the main techniques in the proof are exponential
approximation and contraction principle (\cite{dz92}).

Let $U=U(\theta)$ be a $Beta(1,\theta)$ random variable, $E=[0,1]$, and $\lambda(\theta)=(-\log(\theta))^{-1}$.

\begin{lemma}\label{l2.1}
The family of the laws of $U(\theta)$ satisfies a large deviation principle on $E$
with speed $\lambda(\theta)$ and rate function
\be\label{ratef1}
      I(p)= \left\{\begin{array}{ll}
      0,& p =1\\
      1,& \mb{else.}
      \end{array}\right.
      \ee
\end{lemma}
\proof For any $a<b$ in $E$, let ${\bf I}$ denote one of the intervals $(a,b), [a,b), (a,b],$ and $[a,b]$. It follows from
direct calculation that for $b <1$
\[
\lim_{\theta \ra 0}\lambda(\theta)\log P\{U \in {\bf I}\}
=-\lim_{\theta \ra 0}\frac{\log(1-c^{\theta})}{\log(\theta)}=-1,
\]
where $c=\frac{1-b}{1-a}$. If $b=1$, then $\lim_{\theta \ra 0}\lambda(\theta)\log P\{U \in J\}=0$. These, combined with
compactness of $E$, implies the result.

\hfill $\Box$

Next let $U_1,U_2,\ldots$ be i.i.d. copies of $U$ and
\be\label{gem}
 X_1=U_1,\  X_m=(1-U_1)\cdots(1-U_{m-1})U_m, \ m \geq 2.
\ee

      \begin{lemma}\label{l2.2}
      For any $n \geq 1$, the family of the laws of $P_{1,n}(\theta) = \max\{X_1, \cdots, X_n\}$
      satisfies a large deviation principle on $E$ with speed $\lambda(\theta)$ and rate function
      \be\label{ratef2}
      I_n(p)= \left\{\begin{array}{ll}
      0,& p =1\\
      k,& p \in [\frac{1}{k+1},\frac{1}{k}), k=1,2,\ldots,n-1\\
      n,& \mb{else.}
      \end{array}\right.
      \ee
      \end{lemma}
      \proof Noting that $P_{1,n}(\theta)$ is a continuous function of $(U_1,\ldots,U_n)$, it follows from
      Lemma~\ref{l2.1}, the independency, and the contraction principle that the family of the laws of
      $P_{1,n}(\theta)$ satisfies a large deviation principle on
      $E$ with speed $\lambda(\theta)$ and rate function
      \[
      I'(p)=\inf\{\sum_{i=1}^nI(u_i): u_i \in E, 1\leq i \leq n;
       \max\{u_1, (1-u_1)u_2,\ldots, (1-u_1)\cdots(1-u_{n-1})u_n\}=p\}.
      \]

      For $p=1$, one has $I'(1)=0$ by choosing $u_i=1$ for $i=1,\ldots,n$. If $p$ is in $[1/2,1)$,
then at least one of the
      $u_i$ is not one.
    By choosing $u_1 =p, u_i =1, i=2,\ldots,n$, it follows that  $I'(p)=1$ for $p$ in $[1/2,1)$.

    For each $m\geq 2$, we have
      \beqn
      &&\max\{u_1, (1-u_1)u_2, \ldots, (1-u_1)\cdots(1-u_{m})\}\label{max}\\
      &&=\max\{u_1,
      (1-u_1)\max\{u_2, \ldots, (1-u_2)\cdots(1-u_{m})\}\}.\nn
      \eeqn

      Noting that
      \[
      \max\{u_1,1-u_1\}\geq \frac{1}{2},\  u_1 \in E,
      \]
      it follows from \rf{max} and induction that
      \be\label{max1}
       \max\{u_1, (1-u_1)u_2, \ldots, (1-u_1)\cdots(1-u_{m})\}\geq \frac{1}{m+1}, \  u_i \in E, i=1,\ldots,m.
      \ee

      Thus, for $2\leq k \leq n-1$, and $p$ in $[\frac{1}{k+1},\frac{1}{k})$, in order for the equality
      \[
      \max\{u_1, (1-u_1)u_2,\ldots, (1-u_1)\cdots(1-u_{n-1})u_n\}=p
      \]
      to hold, it is necessary that $u_1,u_2,\ldots, u_k$ are all less than one. In other words,
      $I'(p) \geq k$. Since the function $\max\{u_1, (1-u_1)u_2,\ldots, (1-u_1)\cdots(1-u_{k})\}$ is a surjection
      from $E^k$ into $[\frac{1}{k+1},1]$, there exists $u_1<1,...,u_k<1$ such that
\[
\max\{u_1, (1-u_1)u_2,\ldots, (1-u_1)\cdots(1-u_{k})\}=p.
\]
 By choosing $u_j=1$ for $j=k+1,\ldots,n$, it follows that $I'(p)=k$.

      Finally for $p$ in $[0,\frac{1}{n})$, in order for
      $$\max\{u_1, (1-u_1)u_2,\ldots, (1-u_1)\cdots(1-u_{n-1})u_n\}=p$$
      to have solutions, each $u_i$ has to be less than one and, thus, $I'(p)=n$.
      Therefore, $I'(p)=I_n(p)$ for all $p$ in $E$.

      \hfill $\Box$

      \begin{lemma}\label{l2.3}
      The laws of $P_1(\theta)$ under $PD(\theta)$ satisfy a large deviation principle on $E$ with speed $\lambda(\theta)$
      and rate function
      \be\label{ratef3}
      S_1(p)= \left\{\begin{array}{ll}
      0,& p =1\\
      k,& p \in [\frac{1}{k+1},\frac{1}{k}), k=1,2,\ldots\\
      \infty,& p=0
      \end{array}\right.
      \ee
      \end{lemma}
      \proof Since the law of $(X_1,X_2,\cdots)$ is the same as the size-biased permutation of $PD(\theta)$, it follows
      that $P_1(\theta)$ under $PD(\theta)$ has the same distribution as $\tilde{P}_1(\theta)=\max\{X_i: i\geq 1\}$.
      For any $\delta >0$, it follows from direct calculation that for any $n \geq 1$
      \beq
     P\{\tilde{P}_1(\theta)- P_{1,n}(\theta)>\delta\} &\leq& P\{(1-U_1)\cdots(1-U_n) >\delta\}\\
     &\leq& \delta^{-1} (\frac{\theta}{1+\theta})^n,
      \eeq
     which implies that
     \be\label{eapp1}
    \limsup_{\theta \ra 0}\lambda(\theta)\log P\{\tilde{P}_1(\theta)- P_{1,n}(\theta)>\delta\} \leq -n.
     \ee

     Hence $\{P_{1,n}(\theta):\theta >0\}$ are exponentially good approximations of  $\{\tilde{P}_1(\theta):\theta >0\}$.
    By direct calculation, for every closed subset $F$ of $E$
\[
\inf_{q \in F}S_1(q) =\limsup_{n \ra \infty} \inf_{q\in F} I_n(q).
\]

     This, combined with theorem 4.2.16 in \cite{dz92} and the fact that $S_1(p)$ is a good rate function,
     implies that a large deviation principle holds for the laws of $\tilde{P}_1$
    with speed $\lambda(\theta)$
      and rate function
      \[
       \sup_{\delta>0}\liminf_{n \ra \infty}\inf_{|q-p|<\delta}I_n(q),
      \]
which is clearly equal to $S_1(p)$.

\hfill $\Box$

      For any $m \geq 1$, let
      \be \label{dafe10}
      \nabla_m = \{(p_1,...,p_m):0\leq p_m\leq ...\leq p_1, \sum_{k=1}^m
      p_k \leq 1\},
      \ee
      and set $ Q_{m,\theta}$ to be the law of $(P_1(\theta),...,P_m(\theta))$ under $PD(\theta)$ on space
      $\nabla_m$.

       For any $\delta >0,$ and any $(p_1,...,p_m) \in
      \nabla_m$,
      let
      \beq
      G((p_1,...,p_m);\delta) &=&\{(q_1,...,q_m)\in \nabla_m: |q_k-p_k|<\delta,
      k=1,...,m\},\\
      F((p_1,...,p_m);\delta) &=&\{(q_1,...,q_m)\in \nabla_m:
      |q_k-p_k|\leq\delta, k=1,...,m\}. \eeq

      \begin{lemma}\label{l2.4}
      For fixed $m \geq 2$, the family $\{Q_{m,\theta}:\theta >0\}$
      satisfies a large deviation principle on space $\nabla_m$ with speed $\lambda(\theta)$ and rate function
      \be\label{ratef4}
      S_m(p_1,...,p_m)=
      \left\{\begin{array}{ll}
       0,& (p_1,p_2,...,p_m)=(1,0...,0)\\
       l-1,& 2\leq l \leq m, \sum_{k=1}^l p_k =1, p_l >0 \\
       m + S_1(\frac{p_m}{1-\sum_{i=1}^m p_i}\wedge 1),&  \sum_{k=1}^m p_k <1, p_m>0\\
       \infty,& \mb{else}.
      \end{array}\right.
      \ee
      \end{lemma}
      \proof Let $m \geq 2$ be fixed, and
      $g^{\theta}_1$ denotes the density function of $P_1(\theta)$. Then
      for any $p\in (0,1)$
      \be \label{dafe12}
      g_1^{\theta}(p)p(1-p)^{1-\theta} =\theta \int_0^{(p/(1-p))\wedge
      1}g_1^{\theta}(x)d\,x,
      \ee

      The joint density function $ g^{\theta}_m$ of
      $(P_1(\theta),...,P_m(\theta))$
      is given by (cf.\cite{Wa76})
      \[
      g^{\theta}_m(p_1,...,p_m)=\frac{\theta^{m-1}(1-\sum_{k=1}^{m-1}p_k)^{\theta-
      2}}{p_1\cdots
      p_{m-1}}g^{\theta}_1(\frac{p_m}{1-\sum_{k=1}^{m-1}p_k}),
      \]
      for $$(p_1,...,p_m)\in \nabla_m^{\circ}=\{(p_1,...,p_m)\in
      \nabla_m: 0<p_m<\cdots<p_1<1, \sum_{k=1}^m p_k <1 \},$$ and is zero
      otherwise. Thus for any fixed $(p_1,...,p_m) \in
      \nabla_m^{\circ}$ we have \be\label{den1}
      g_m^{\theta}(p_1,...,p_m)=\frac{\theta^m (1-\sum_{k=1}^m
      p_k)^{\theta-1}}{p_1\cdots p_m}\int_0^{(p_m/(1-\sum_{k=1}^m
      p_k))\wedge 1}g_1^{\theta}(u)d\,u. \ee

\vspace{0.26cm}

The key step in the proof is to show that for every $(p_1,...,p_m)$ in $\nabla_m$,

\beqn
&&\lim_{\delta \ra 0}\liminf_{\theta \ra 0}\lambda(\theta)\log Q_{m,\theta}(F((p_1,...,p_m);\delta))\nn\\
&&\ \ \ \ =\lim_{\delta \ra 0}\limsup_{\theta \ra 0}\lambda(\theta)\log Q_{m,\theta}(G((p_1,...,p_m);\delta))\label{key}\\
&&\ \ \ \ =-S_m(p_1,...,p_m).\nn
\eeqn

For any $(p_1,...,p_m)$ in $\nabla_m$ satisfying $\sum_{i=1}^m p_i >0$, define
      \be\label{index}
       r=r(p_1,\ldots, p_m)= \max\{i: 1\leq i \leq m, p_i >0\}.
      \ee

We divide the proof into several disjoint cases.
\vspace{0.3cm}

      \noindent {\bf Case I:}\ \  $r=1$, i.e., $(p_1,...,p_m)=(1,...,0)$.
\vspace{0.3cm}

       For any $\delta >0$, \[
 F((1,...,0); \delta) \subset \{(q_1,...,q_m)\in \nabla_m: |q_1-1|\leq \delta\},
\]
and one can choose $\delta'< \delta$ such that
\[
\{(q_1,...,q_m)\in \nabla_m: |q_1-1|< \delta'\} \subset G((1,...,0); \delta).
\]

These combined with Lemma~\ref{l2.3} implies \rf{key} in this case.

\vspace{0.3cm}

      \noindent {\bf Case II:}\ \  $r=m, \sum_{k=1}^m p_k <1$.
\vspace{0.3cm}

      Choose $\delta >0$ so that
       $$ \delta < \min\{p_m,\frac{1-\sum_{i=1}^m p_i}{m}\}.$$

        By \rf{den1},
      we have that for any $(q_1,...,q_m)$ in $F((p_1,...,p_m), \delta)\cap \nabla_m^{\circ}$
      \[
      g^{\theta}_m(q_1,...,q_m)\leq\frac{\theta^{m}(1-\sum_{k=1}^{m}(p_k+\delta))^
      {\theta-1}}{(p_1-\delta)\cdots
      (p_{m}-\delta)}\int_0^{\frac{p_m+\delta}{1-\sum_{k=1}^m
      (p_k+\delta)}\wedge 1}g_1^{\theta}(u)d\,u ,
      \]
      which, combined with Lemma~\ref{l2.3}, implies
      \beqn
      &&\lim_{\delta \ra 0}\limsup_{\theta \ra
      0}\lambda(\theta)\log
      Q_{m,\theta}\{F((p_1,...,p_m);\delta)\}\nn\\
      && \ \ \ \ \leq
      -m +\lim_{\delta \ra 0}\limsup_{\theta \ra
      0}\lambda(\theta)\log
      P\{P_1(\theta)\leq \frac{p_m+\delta}{1-\sum_{k=1}^m
      (p_k+\delta)}\wedge 1\}\label{ub3}\\
      && \ \ \ \ \leq -[m +S_1(\frac{p_m}{1-\sum_{i=1}^m p_i}\wedge 1)]\nn,
      \eeqn
      where the right continuity of $S_1(\cdot)$ is used in the last inequality.

       On the other hand, let
       \[
       \tilde{G}(p_1,...,p_m), \delta)=\prod_{i=1}^m(p_i +\frac{\delta}{2},p_i+\delta)\cap \nabla_m^{\circ},
       \]
       which is clearly a subset of $G((p_1,...,p_m), \delta)$.
       Using \rf{den1} again it follows that for any
      $(q_1,...,q_m)$ in $\tilde{G}((p_1,...,p_m), \delta)$
      \[
      g^{\theta}_m(q_1,...,q_m)\geq \theta^{m}\frac{(1-\sum_{k=1}^{m}(p_k
      +\delta/2))^{\theta-1}}{(p_1+\delta)\cdots(p_m+\delta)} \int_0^{((p_m+\delta/2)/(1-\sum_{k=1}^m
      (p_k+\delta/2)))\wedge 1}g_1^{\theta}(u)d\,u,
      \]
      which, combined with Lemma~\ref{l2.3}, implies
      \beq
      \liminf_{\theta \ra 0}\lambda(\theta)\log
      Q_{m,\theta}\{G((p_1,...,p_m);\delta)\}
      &\geq& \liminf_{\theta \ra 0}\lambda(\theta)\log
      Q_{m,\theta}\{\tilde{G}((p_1,...,p_m);\delta)\}\\
      &\geq& -m -S_1(\frac{p_m+\delta/2}{1-\sum_{i=1}^m (p_i+\delta/2)} \wedge 1).
      \eeq
      It follows, by letting $\delta$ go to zero, that
        \be\label{lb4}
      \liminf_{\delta \ra 0}\liminf_{\theta\ra
      \infty}\frac{1}{\theta}\log
      Q_{m,\theta}\{G((p_1,...,p_m);\delta)\}\geq -S_m(p_1,...,p_m).
      \ee

     \vspace{0.3cm}

      \noindent {\bf Case III:}\ \  $2\leq r \leq m-1, \sum_{i=1}^r p_i<1$ or $p_1=0$.

      \vspace{0.3cm}

 This case follows from estimate \rf{ub3} and the fact that $S_1(0)=-\infty$.

\vspace{0.3cm}

      \noindent {\bf Case IV:}\ \  $r=m, \sum_{k=1}^m p_k =1$.
\vspace{0.3cm}

  Noting that for any $\delta>0$
  \[
  F((p_1,...,p_m);\delta)\cap\nabla^{\circ}_m \subset \{(q_1,...,q_m)\in \nabla_m^{\circ}:
   |q_i-p_i|\leq \delta, i=1,...,m-1\}.
  \]

 By applying {\bf Case II} to
 $(P_1(\theta),...,P_{m-1}(\theta))$ at the point $(p_1,...,p_{m-1})$, we get
      \be\label{ub4}
\lim_{\delta \ra 0}\limsup_{\theta \ra
      0}\lambda(\theta)\log
      Q_{m,\theta}\{F((p_1,...,p_m);\delta)\} \leq -[m-1 +S_1(1)]=-(m-1).
      \ee

 On the other hand, one can choose $\delta>0$ small so that
$
\frac{q_m}{1-\sum_{i=1}^m q_i} >1
$
for any $(q_1,...,q_m)$ in $
 G((p_1,...,p_m);\delta)\cap \nabla_m^{\circ}$.

Set
\[
\tilde{G}=\{(q_1,...,q_m)\in \nabla_m^{\circ}: p_i<q_i< p_i +\delta/(m-1), i=1,...,m-1; p_m -\delta < q_m <p_m\}.
\]

Clearly $\tilde{G}$ is a subset of $G((p_1,...,p_m);\delta)$. It follows from \rf{den1} that
for any $(q_1,...,q_m)$ in $\tilde{G}$,

\[
g_m^{\theta}(q_1,...,q_m)\geq \frac{\theta^{m-1}[\theta(1-\sum_{i=1}^m q_i)^{\theta-1}]}
{(p_1+\delta/(m-1))\cdots (p_{m-1}+\delta/(m-1))p_m}.
\]

For $m \geq 2$, let
\[
A_m =\{(q_1,...,q_{m-1})\in \nabla_{m-1}:p_i<q_i< p_i +\delta/(m-1), i=1,...,m-1, \sum_{j=1}^{m-1} q_j <1\}
\]
Then
\beq
&&\int_{\tilde{G}}\theta(1-\sum_{i=1}^m q_i)^{\theta-1}dq_1\cdots dq_m\\
&& \hspace{1.5cm}=\int_{A_m}dq_1\cdots dq_{m-1} \int_{p_m-\delta}^{p_m \wedge (1-\sum_{i=1}^{m-1}q_i)}
\theta(1-\sum_{i=1}^m q_i)^{\theta-1}dq_m\\
&& \hspace{1.5cm}=\int_{A_m}(1+\delta-p_m-\sum_{i=1}^{m-1}q_i)^{\theta} dq_1\cdots dq_{m-1},
\eeq
which converges to a strictly positive number depending only on $\delta$ and $(p_1,...,p_m)$  as $\theta$ goes to zero.
Hence
\be\label{newlb4}
\lim_{\delta \ra 0}\liminf_{\theta \ra
      0}\lambda(\theta)\log
      Q_{m,\theta}\{G((p_1,...,p_m);\delta)\}\geq  \lim_{\delta \ra 0}\liminf_{\theta \ra
      0}\lambda(\theta)\log
      Q_{m,\theta}\{\tilde{G}\}\geq -(m-1).
\ee

\vspace{0.3cm}

\noindent {\bf Case V:} $2\leq r\leq m-1, \sum_{i=1}^r p_i=1$.

\vspace{0.3cm}

First note that for any $\delta >0$, $ F((p_1,...,p_m);\delta)$ is a subset of
$$ \{(q_1,...,q_m)\in \nabla_m: |q_i-p_i|\leq \delta, i=1,...,r\}.$$

On the other hand, for each $\delta >0$ one can choose $\delta_0< \delta$ such that for any $\delta'\leq \delta_0$
\[
G((p_1,...,p_m);\delta) \supset \{(q_1,...,q_m)\in \nabla^{\circ}_m; |q_i-p_i|<\delta',i=1,...,r\}.
\]

Thus the result now follows from {\bf Case IV} for $(P_1(\theta),...,P_r(\theta))$.

The lemma now follows from \rf{key} and the fact that
        $\nabla_m$ is compact.

      \hfill $\Box$

      For any $n \geq 1$, set
      \[
      L_n = \{(p_1,...,p_n,0,0,...) \in \nabla: \sum_{i=1}^n p_i = 1\}
      \]
      and
      \[
        L= \bigcup_{i=1}^{\infty}L_i.
      \]

      Now we are ready to state and prove the main result of this section.

      \begin{theorem}\label{P3}
      The family $\{PD(\theta):\theta >0\}$ satisfies a large deviation principle with speed
      $\lambda(\theta)$ and rate function

      \be\label{rate function3}
      S({\bf p})=
      \left\{\begin{array}{ll}
      0, & {\bf p}\in L_1\\
      n-1,& {\bf p} \in L_n, p_n >0, n\geq 2\\
      \infty,& {\bf p}\not\in L
      \end{array}\right.
      \ee
      \end{theorem}
      \proof First note that the topology of the space
      $\nabla$ can be generated by the following metric
      \[
      d({\bf p}, {\bf q})=\sum_{k=1}^{\infty}\frac{|p_k-q_k|}{2^k},
      \]
      where ${\bf p}=(p_1,p_2,...), {\bf q}=(q_1,q_2,...)$. For any fixed $\delta >0$, let $B({\bf
      p},\delta)$ and $\bar{B}({\bf p},\delta)$ denote the respective
      open and closed balls centered at ${\bf p}$ with radius $\delta
      >0$.
\vspace{0.3cm}

      We start with the case that ${\bf p}$ is not in $L$.

  \vspace{0.3cm}

      For any $k \geq 1, \delta'>0$, set

      \[
      \bar{B}_{k,\delta'}({\bf p})=\{(q_1,q_2,...) \in \nabla:
      |q_i-p_i|\leq \delta', i=1,...,k\}.
      \]

      Choose $\delta >0$ so that $2^k \delta <\delta'$. Then
      \[
      \bar{B}({\bf p},\delta)\subset \bar{B}_{k,\delta'}({\bf p}),
      \]
      and
      \beqn
      \lim_{\delta \ra 0}\limsup_{\theta \ra 0}\lambda(\theta)\log
       PD(\theta)\{\bar{B}({\bf p},\delta)\} &\leq &  \limsup_{\theta \ra 0}\lambda(\theta)\log
       PD(\theta)\{\bar{B}_{k,\delta'}({\bf p})\}\nn\\
       &\leq& \limsup_{\theta \ra 0}
       \lambda(\theta)\log Q_{k,\theta}\{F((p_1,...,p_k),
      \delta')\}\label{ub1}\\
      & \leq &-\inf\{S_k(q_1,...,q_k):(q_1,...,q_k)\in F((p_1,...,p_k),
      \delta')\}.\nn \eeqn
      Letting $\delta'$ go to zero, and then $k$
      go to infinity, we get
      \be\label{ub2}
      \lim_{\delta \ra
      0}\liminf_{\theta \ra 0}\lambda(\theta)\log
      PD(\theta)\{B({\bf p},\delta)\} =\lim_{\delta \ra
      0}\limsup_{\theta \ra 0}\lambda(\theta)\log
      PD(\theta)\{\bar{B}({\bf p},\delta)\} =-\infty. \ee

\vspace{0.3cm}

      Next consider the case of ${\bf p}$ belonging to $L$. Without loss of generality, we assume that
      ${\bf p}$ belongs to $L_n$ with $p_n >0$.

\vspace{0.3cm}

       For any $\delta >0,$
      let
      \beq
      \tilde{G}({\bf p};\delta) &=&\{{\bf q}\in \nabla: |q_k-p_k|<\delta,
      k=1,...,n\},\\
      \tilde{F}({\bf p};\delta) &=&\{{\bf q}\in \nabla:
      |q_k-p_k|\leq\delta, k=1,...,n\}. \eeq

      Clearly, $\bar{B}({\bf p},\delta)$ is a subset of $ \tilde{F}({\bf p};2^n\delta)$.
Since $\sum_{i=1}^n p_i =1$, it follows that, for any $\delta >0$, one can find
$\delta' <\delta$ such that
\[
B({\bf p},\delta)\supset \tilde{G}({\bf p};\delta').
\]

Using results on $(P_1(\theta),...,P_n(\theta))$ in {\bf Case V} in the proof of
Lemma~\ref{l2.4}, we get
\be\label{key1}
\lim_{\delta \ra 0}\liminf_{\theta \ra 0}\lambda(\theta)\log PD(\theta)(B({\bf p}, \delta))
=\lim_{\delta \ra 0}\limsup_{\theta \ra 0}\lambda(\theta)\log PD(\theta)(\bar{B}({\bf p},\delta))
=-(n-1).
\ee

Finally, the theorem follows from the compactness of $\nabla$.

      \hfill $\Box$

      \vspace{0.4cm}

      \noindent
      {\bf Remarks.} {\rm 1.} Consider the rate function $S(\cdot)$ as an ``energy" function, then the energy
      needed to get $n \geq 2$ different alleles is $n-1$. The values of $S(\cdot)$ form a``ladder of energy".
      The energy needed to get infinite number of alleles is infinity and thus it is impossible to have
      infinitely many alleles under large deviation.

      {\rm 2.} The effective domain of $S(\cdot)$, defined as  $\{{\bf p} \in \nabla:
      S({\bf p})<\infty\}$, is clearly $L$.  This is in sharp contrast to the result in
 \cite{dawson-feng06} where the rate function associated with large mutation rate has an effective
 domain of $\{{\bf p} \in \nabla: \sum_{i=1}^{\infty}p_i <1\}$. The two effective domains are disjoint.
 One is part of the boundary of $\nabla$ and the other is the interior of $\nabla$, and both have no
 intersections with the set
$\{{\bf p}\in \nabla: p_1>p_2\cdots>0, \ \sum_{i=1}^{\infty}p_i=1\}$.

      {\rm 3.} The large deviation result with small mutation rate is limited to $PD(\theta)$.
      Both the Dirichlet process and the two-parameter Poisson-Dirichlet distribution converge to
      random limits when mutation rate becomes small. Thus
      it is not clear how to formulate large deviation problems for them when the mutation rate approaches zero.

      \section{Applications}

       In this section we will discuss two applications of Theorem~\ref{P3}. The first one is concerned with the
       large deviation principle
      for the homozygosity.

      A random sample of size $r\geq 2$ is selected from a population
      whose allelic types have distribution $PD(\theta)$. The
      probability that all samples are of the same type is called the
      $r$th order population homozygosity and is given by
      \be \label{dafe20}
      H_r(P_1(\theta),...)=\sum_{i=1}^{\infty}P_i^r(\theta).
      \ee

      It is clear that $H_r(P_1(\theta),...)$ converges to one as $\theta$ approaches zero. Our next theorem describes the
      large deviations of $H_r(\theta)$ from one.

      \begin{theorem}\label{applt1}
      For any integer $r\geq 2$, the family of the laws of $H_r(P_1(\theta),...)$ satisfies a large deviation principle on $E$ with speed
      $\lambda(\theta)$ and rate function

      \be\label{ratef5}
      J(p)= \left\{\begin{array}{ll}
      0,& p =1\\
      n-1,& p \in [\frac{1}{n^{r-1}},\frac{1}{(n-1)^{r-1}}), n=2,\ldots\\
      \infty,& p=0
      \end{array}\right.
      \ee
      Thus in terms of large deviations, $H_r(P_1(\theta),...)$  behaves the same as $P^{r-1}_1(\theta)$.
      \end{theorem}
      \proof For any integer $r>1$, $H_r({\bf p})$ is clearly continuous on $\nabla$.
            By Theorem~\ref{P3} and the
      contraction principle, the family of the laws of $H_r(P_1(\theta),...)$
      satisfies a large deviation principle with speed $\lambda(\theta)$ and rate function
      \[
      \inf\{S({\bf q}):{\bf q} \in \nabla, H_r({\bf
      q})=p\}= \inf\{S({\bf q}):{\bf q} \in L, H_r({\bf
      q})=p\}.
      \]

      For $p=1$, it follows by choosing ${\bf q}=(1,0,...)$ that $\inf\{S({\bf q}):{\bf q} \in \nabla, H_r({\bf
      q})=p\}=0$. For $p=0$, there does not exist ${\bf q}$ in $L$ such that $H_r({\bf
      q})=p$. Hence $\inf\{S({\bf q}):{\bf q} \in L, H_r({\bf
      q})=p\}=\infty.$

      For any $n \geq 2$, the minimum of $\sum_{i=1}^n q^r_i$ over $L_n$ is $n^{-(r-1)}$
      which is achieved when all $q_i's$ are equal. Hence for

\[
\ p \in [n^{-(r-1)}, (n-1)^{-(r-1)}),
\]
we have
\[
\inf\{S({\bf q}):{\bf q} \in \nabla, H_r({\bf
      q})=p\}= n-1 =J(p).
\]

 \hfill       $\Box$

 For any
  $\alpha(\theta) >0$ and any nonzero constant $s$, the Poisson-Dirichlet distribution with symmetric selection is
  a probability measure on $\nabla$
   given by
\[
P_{\alpha(\theta), s H_r, \theta}(d{\bf p})=(\int_{\nabla}e^{s \alpha(\theta)H_r({\bf q})}PD(\theta)(d{\bf q}))^{-1}
e^{s \alpha(\theta)H_r({\bf p})}PD(\theta)(d{\bf p}),
\]
where $\alpha(\theta)$ is the selection intensity, $s>0 (<0)$ corresponds to underdominant (overdominant) selection.

In our second application, Theorem~\ref{P3} is used to derive the large deviation principle for
the Poisson-Dirichlet distribution with symmetric selection.

\begin{theorem}\label{applt2}
      The family $\{P_{\alpha(\theta), s H_r, \theta}: \theta >0\}$ satisfies a large deviation principle on
      $\nabla$ with speed $\lambda(\theta)$ and rate function

      \be\label{ratef6}
      S'({\bf p})= \left\{\begin{array}{ll}
      S({\bf p}),& \lim_{\theta\ra 0}\alpha(\theta)\lambda(\theta)=0\\
      S({\bf p})+ s(1-H_r({\bf p})),& \alpha(\theta)\lambda(\theta)=1, s>0\\
      S({\bf p})+|s|H_r({\bf p})-\inf\{\frac{|s|}{n^{r-1}}+n-1: n\geq 1\},& \alpha(\theta)\lambda(\theta)=1, s<0
      \end{array}\right.
      \ee
      \end{theorem}

\proof Theorem~\ref{P3}, combined with Varadhan's lemma and the Laplace method, implies that the
family $\{P_{\alpha(\theta), s H_r, \theta}: \theta >0\}$ satisfies a large deviation principle on
      $\nabla$ with speed $\lambda(\theta)$ and rate function
\[
\sup\{sH_r({\bf q})-S({\bf q}): {\bf q} \in \nabla\}- (sH_r({\bf p})-S({\bf p})).
\]

The theorem then follows from the fact that
\[
\sup\{sH_r({\bf q})-S({\bf q}): {\bf q} \in \nabla\}= \left\{\begin{array}{ll}
      s,&  s>0\\
      -\inf\{\frac{|s|}{n^{r-1}}+n-1: n\geq 1\},&  s<0
      \end{array}\right.
\]

\hfill $\Box$

\noindent {\bf Remark.} It is clear from the above theorem that selection has an impact on large deviations
only when the selection intensity is comparable to $\lambda(\theta)^{-1}$. Assuming that
$\alpha(\theta)=(\lambda(\theta))^{-1}$. When $s>0$, homozygote has selective advantage. Thus the small mutation rate limit is
$(1,...)$ and large deviations from it become more difficult comparing to the neutral case. This is confirmed
through the fact that $S'({\bf p})$ is bigger than $S({\bf p})$. For $s<0$, heterozygote has selection advantage.
The fact that $S'(\cdot)$ may reach zero at a point that is different from $(1,...)$ shows that several alleles can coexist
in the population when the overdominant selection goes to infinity and the mutation approaches zero. In some cases
such as $r=2, s=-k(k+1), k\geq 1$,
$S'(\cdot)$ is zero at more than one point.

\end{document}